\documentclass[twoside]{article}
\title{Gorenstein Global Dimensions and Cotorsion Dimension of Rings}
\date{}

\usepackage{amsfonts}
\usepackage{amsmath}
\usepackage{amssymb}
\usepackage{latexsym}
\usepackage{graphicx}

\newtheorem{thm}{\bf Theorem}[section]

\newtheorem{lem}[thm]{\bf Lemma}
\newtheorem{prop}[thm]{\bf Proposition}
\newtheorem{defn}[thm]{\bf Definition}

\newtheorem{exmp}[thm]{\bf Example}

\catcode`\ç=13
\defç{\c{c}}
\catcode`\é=13
\defé{\'e}
\catcode`\à=13
\defà{\`a}
\catcode`\è=13
\defè{\`e}
\catcode`\â=13
\defâ{\^a}
\catcode`\ù=13
\defù{\`u}
\catcode`\ê=13
\defê{\^e}
\catcode`\î=13
\defî{\^\i}
\catcode`\ô=13
\defô{\^o}

\def\proof{{\parindent0pt {\bf Proof.\ }}}

\def\wdim{{\rm wdim}}
\def\gldim{{\rm gldim}}

\def\Gwdim{{\rm G\!-\!wdim}}
\def\Ggldim{{\rm G\!-\!gldim}}
\def\cot.D{{\rm cot.D}}

\def\FPD{{\rm FPD}}

\def\pd{{\rm pd}}
\def\fd{{\rm fd}}
\def\id{{\rm id}}
\def\FPid{{\rm FP\!-\!id}}
\def\cd{{\rm cd}}

\def\g{{\rm G\!-\!dim}}
\def\Gpd{{\rm Gpd}}
\def\Gfd{{\rm Gfd}}
\def\Gid{{\rm Gid}}

\def\Ext{{\rm Ext}}
\def\Tor{{\rm Tor}}
\def\Hom{{\rm Hom}}

\def\sup{{\rm sup}}

\newcommand{\cqfd}
{\hspace{1cm}
\rule{2mm}{2mm}%
\medbreak%
\par%
}

\begin{document}
\thispagestyle{empty}
\maketitle \vspace*{-1.5cm}
\begin{center}{\large\bf Driss Bennis and Najib Mahdou}

\bigskip

 \small{Department of Mathematics, Faculty of Science and Technology of Fez,\\ Box 2202, University S. M.
Ben Abdellah Fez, Morocco\\[0.12cm]
driss\_bennis@hotmail.com\\
mahdou@hotmail.com}
\end{center}

\bigskip\bigskip
\noindent{\large\bf Abstract.} In this paper,  we establish, as a
generalization of a result on the classical homological dimensions
of commutative rings, an upper bound on the Gorenstein global
dimension of commutative rings using the global cotorsion
dimension of rings. We use this result to compute the Gorenstein
global dimension of some particular cases of trivial extensions of
rings and of group rings.\bigskip

\bigskip\bigskip



\begin{section}{Introduction} Throughout this paper    all   rings
are commutative with identity element and all modules are
unitary.\\
\indent For a ring $R$  and   an $R$-module $M$, we use
$\pd_R(M),\ \id_R(M)$, and $\fd_R(M)$ to denote, respectively, the
classical projective, injective and flat dimensions of $M$. By
$\gldim(R)$ and $\wdim(R)$ we  denote, respectively, the classical
global and weak global dimensions of $R$.\bigskip

The Gorenstein homological dimensions theory originated in the
works of Auslander and Bridger  \cite{A1} and \cite{A2}, where
they introduced the G-dimension, $\g_R(M)$, of any finitely
generated module $M$ and over any Noetherian ring $R$. The
G-dimension is analogous to the classical projective dimension and
shares some of its principal properties (see \cite{LW} for more
details). However, to complete the analogy an extension of the
G-dimension to non-necessarily finitely generated modules is
needed. This is done in \cite{GoIn, GoInPj}, where the Gorenstein
projective dimension was defined over arbitrary rings (as an
extension of the G-dimension to modules that are not necessarily
finitely generated), and the Gorenstein injective dimension was
defined as a dual notion of the Gorenstein projective dimension.
And also to complete the analogy with the classical homological
dimensions theory, the Gorenstein flat dimension was introduced in
\cite{GoPlat}. Since then, several results on the classical
homological dimensions were extended to the Gorenstein homological
dimensions. Namely, the majority of works on the Gorenstein
homological dimensions attempt to confirm the following
meta-theorem (please see Holm's thesis \cite[page v]{HPhD}):
``\textit{Every result in classical homological algebra has a
counter part in Gorenstein homological algebra.}" (for more
details see also \cite{LW, CFH, Rel-hom, HH}). In line with this,
the Gorenstein global dimensions of commutative rings were
investigated in \cite{BM1} (and \cite{BM2}). It is proved, for any
ring $R$ \cite[Theorems 3.1.3 and 3.2.1]{BM1}:
$$\sup\{\Gfd_R(M)\,|\,M\;R\!-\!module\}\leq
\sup\{\Gpd_R(M)\,|\,M\;R\!-\!module\}= \sup
\{\Gid_R(M)\,|\,M\;R\!-\!module\}.$$  So, according to the
terminology of the classical theory of homological dimensions of
rings,  the common value of
$\sup\{\Gpd_R(M)\,|\,M\;R\!-\!module\}$ and $\sup
\{\Gid_R(M)\,|\,M\;R\!-\!module\}$ is called \textit{Gorenstein
global dimension} of $R$, and denoted by $\Ggldim(R)$, and the
homological invariant  $\sup\{\Gfd_R(M)\,|\,M\;R\!-\!module\}$ is
called \textit{Gorenstein weak global dimension} of $R$, and
denoted by $\Gwdim(R)$.\\
\indent The Gorenstein weak global and global dimensions are
refinements of the classical  weak and global dimensions of rings,
respectively; that is \cite[Propositions 3.11 and 4.5]{BM1}:
$\Ggldim(R)\leq \gldim(R)$ and $\Gwdim(R)\leq \wdim(R)$, with each
of the two inequalities becomes equality if
$\wdim(R)$ is finite.\\
\indent If $R$ is a Noetherian ring, then \cite[Corollary
2.3]{BM1}: $\Gwdim(R)=\Ggldim(R)$, such that:
$$\begin{array}{l}
  \Ggldim(R)\leq n \quad \Longleftrightarrow \quad R \ \mathrm{is} \ n\!-\!\mathrm{Gorenstein}.\\
\end{array}
$$
Recall that a ring $R$ is said to be $n$-Gorenstein, for a
positive integer $n$, if it is Noetherian  with self-injective
dimension less  or equal than $n$ (i.e., $\id_R(R)\leq n$); and
$R$ is said to be Iwanaga-Gorenstein, if it is $n$-Gorenstein for
some positive integer $n$ (please see \cite[Section
9.1]{Rel-hom}). Notice that
$0$-Gorenstein rings are the well-known quasi-Frobenius rings.\\
\indent For a coherent ring $R$, we have \cite[Theorem 4.11]{BM1}:
$$\begin{array}{l}
\Gwdim(R)\leq n \quad \Longleftrightarrow \quad R\ \mathrm{is} \ n\!-\!FC.\\
 \end{array}
$$
Recall that a ring $R$ is said to be $n$-FC, for a positive
integer $n$, if it is coherent and $\FPid_R(R)\leq n$
\cite{FCring}; where $\FPid_R(M)$ denotes, for an $R$-module $M$,
the FP-injective dimension, which is defined to be the least
positive integer $n$ for which $\Ext_{R}^{n+1}(P,M)=0$ for all
finitely presented $R$-modules $P$. Notice that the $0$-FC rings
coincide (in commutative setting) with the IF-rings; i.e., rings
over which every injective module is flat (please see \cite{Colby,
Jain, Stenst}).\bigskip

In this paper, we continue the study of the  Gorenstein global
dimensions of commutative rings started in \cite{BM1} and
\cite{BM2}. The paper extends some results on the classical global
homological dimensions to the  Gorenstein global dimensions. To
see that recall the following:\medskip

In  \cite{DM}, Ding and Mao introduced the cotorsion dimension of
modules and rings, which are defined as follows:

\begin{defn}[\cite{DM}]\label{def-cot-dim}Let $R$ be a ring.\\
\indent The cotorsion dimension of an $R$-module $M$, denoted by
$\cd_R(M)$, is the least positive integer $n$
for which $\Ext^{n+1}_{R}(F, M)=0$ for all flat $R$-modules $F$.\\
\indent  The global cotorsion dimension of  $R$, denoted by
$\cot.D(R)$, is defined as the supremum of the cotorsion
dimensions of $R$-modules.
\end{defn}

The global cotorsion dimension of rings measures how far away a
ring is from being perfect: the perfect rings are those rings over
which every flat module is projective (please see \cite{bass}).
Namely, we have, for a ring $R$ and a positive integer $n$,
$\cot.D(R)\leq n$ if and only if every flat $R$-module $F$ has
projective dimension less or equal than $n$ \cite[Theorem 7.2.5
(1)]{DM}. In \cite{n-perfect}, a ring that satisfies the last
condition is called $n$-perfect. So, we have: $\cot.D(R)\leq n$ if
and only if $R$ is $n$-perfect. Particularly, $\cot.D(R)=0$ if and
only if $R$ is
$0$-perfect  if and only  if  $R$ is perfect.\\
\indent The global cotorsion dimension of rings is also used to
give an upper bound on the global dimension of rings as follows
\cite[Theorem 7.2.11]{DM}: For any ring $R$,  we have the
inequality:
$$ \gldim(R)\leq \wdim(R)+ \cot.D(R).\\$$

The main result of this paper (Theorem \ref{thm-main}) extends
this inequality to the  Gorenstein global dimensions of coherent
rings. This result enables us to compute the Gorenstein global
dimension of a particular case of trivial extensions of rings
(Proposition \ref{pro-Ggldim-trivial}).\\
\indent In the end of the paper, we investigate the global
cotorsion dimension of group rings (Theorem \ref{thm-RG-Cot}).
This is used with the main result to compute the Gorenstein global
dimension of a particular case of group rings (Proposition
 \ref{prop-RG-appl}).

\end{section}
\bigskip\bigskip

 \begin{section}{Main results}
Our main result is the following:

\begin{thm}\label{thm-main}
If $R$ is a coherent ring, then:
$$\cot.D(R)\leq\Ggldim(R)\leq\Gwdim(R)+\cot.D(R).$$
\indent In particular:
 \begin{itemize}
    \item If $\cot.D(R)=0$ (i.e., $R$ is perfect), then $\Gwdim(R)= \gldim(R)$.
    \item  If $\Gwdim(R)=0$ (i.e., $R$ is an IF-ring), then $\cot.D(R)=
\Ggldim(R)$.
 \end{itemize}
\end{thm}

To prove this theorem, we need the following result, which is a
generalization of the characterization of the Gorenstein
projective dimension over Iwanaga-Gorenstein rings \cite[Theorem
2.1]{Gflatcover}.

\begin{lem}\label{lem-Gpro-n-perfect-FC}
Let $R$ be both an $n$-FC ring and an $m$-perfect ring, where $n$
and $m$ are positive integers. For any $R$-module $M$,
we have the following equivalence, for a positive integer $k$:\\
$\Gpd_R(M)\leq k \ \Leftrightarrow \ \Ext^{j}_{R}(M,P)=0 $ for all
$j\geq k+1$ and all modules $P$ with finite $\pd_R(P)$.
\end{lem}

In the proof of this lemma we use the notion of a flat preenvelope
of modules which is defined as follows:

\begin{defn}[\cite{Rel-hom}]\label{def-prenvel}
Let $R$ be a ring and let $F$ be a flat $R$-module. For an
$R$-module $M$, an homomorphism $\varphi :\, M\rightarrow F$ is
called  a flat preenvelope, if for any homomorphism $\varphi' :\,
M\rightarrow F'$ with $F'$ is a flat module, there is an
homomorphism $f :\, F\rightarrow F'$ such that $\varphi'
=f\varphi$.
\end{defn}

The coherent rings can be characterized by the notion  of a flat
preenvelope of modules  as follows:

\begin{lem}[\cite{Rel-hom}, Proposition  6.5.1]\label{car-coherent-prevelo}
A ring $R$ is coherent if and only if every $R$-module has a flat
preenvelope.
\end{lem}

\noindent\textbf{Proof of Lemma \ref{lem-Gpro-n-perfect-FC}.} The
direct implication holds over arbitrary rings
by \cite[Theorem 2.20]{HH}.\\
\indent Conversely, consider an exact sequence of $R$-modules:
$$0\rightarrow K_n \rightarrow P_{n-1}\rightarrow\cdots
\rightarrow P_0\rightarrow M \rightarrow 0$$ where each $P_i$ is
projective. We have $\Ext^{n+i}(M,Q)\cong\Ext^{i}(K_n,Q)$ for all
$i\geq 1$ and all modules $ Q$. Then, to prove this implication,
it is sufficient to prove it for $k=0$. Then, we assume that
$\Ext^i(M,P)=0$ for all $i\geq 1$ and all $R$-modules $P$  with
finite projective dimension, and we prove that $M$ is Gorenstein
projective. This is equivalent to prove, from \cite[Proposition
2.3]{HH}, that there exists an exact sequence of $R$-modules:
$$\alpha=\quad 0\rightarrow M\rightarrow P^0 \rightarrow P^1
\rightarrow \cdots,$$ where each $P^i$ is projective, such that
$\Hom_R ( -, P) $ leaves the
sequence $\alpha$ exact whenever $P$ is a projective $R$-module.\\
\indent The proof of this implication is analogous to the one of
\cite[Theorem 2.1 ($1\Rightarrow4$)]{Gflatcover}. For
completeness, we give a proof here.\\
\indent As usual (see for instance the proofs of \cite[Theorems
4.2.6 and 5.1.7]{LW}), to construct the sequence $\alpha$,  it is
sufficient to prove the existence of a short exact sequence of
$R$-modules:
$$0 \rightarrow M\rightarrow P^0 \rightarrow G^0 \rightarrow 0,$$
where $P^0$ is projective, such that $\Ext^i(G^0,P)=0$ for all
$i>0$ and all $R$-modules $P$ with finite projective dimension
(and then the sequence $\alpha$  is recursively constructed).\\
\indent First, we prove that $M$ can be embedded into a flat
$R$-module. For that, pick a short exact sequence of $R$-modules
$0 \rightarrow M\rightarrow I\rightarrow E \rightarrow 0$, where
$I$ is injective. For this $I$ pick a short exact sequence of
$R$-modules $0 \rightarrow Q\rightarrow P\rightarrow I\rightarrow
0$, where $P$ is projective. Consider the following pullback
diagram:
$$\begin{array}{ccccccccc}
      &   & 0 &  & 0 & &  &   &   \\
            &   & \downarrow&   & \downarrow&   &  &   &   \\
     &   & Q& =\!= & Q &  &   &  &   \\
      &   & \downarrow&   & \downarrow&   &  &   &   \\
   0 & \rightarrow & D& \rightarrow & P & \rightarrow & E & \rightarrow & 0 \\
   &   & \downarrow &   & \downarrow &   &||&   &   \\
   0 & \rightarrow & M& \rightarrow & I& \rightarrow & E & \rightarrow & 0 \\
      &   & \downarrow&   & \downarrow&   &  &   &   \\
        &   & 0 &  & 0 & &  &   &   \\
  \end{array}$$
Since $I$ is injective,  $\pd(I)<\infty$ (From \cite[Theorem
4.11]{BM1}  and since $R$ is $m$-perfect). Then, $\pd(Q)<\infty$.
By hypothesis, $\Ext(M,Q)=0$, and then the first vertical exact
sequence is split, so $M$ embeds into $D$ which is an
$R$-submodule of the projective (then flat) $R$-module $P$.\\
\indent  The fact that $M$ embeds  into a flat $R$-module implies,
from Lemma \ref{car-coherent-prevelo} and Definition
\ref{def-prenvel}, that $M$ admits an injective flat preenvelope
$\varphi :\, M \rightarrow F$. For such flat $R$-module $F$,
consider a short exact sequence of $R$-modules $0 \rightarrow H
\rightarrow P^0 \stackrel{f}\rightarrow F \rightarrow 0$, where
$P^0$ is projective, then $H$ is a flat $R$-module, hence it has
finite projective dimension (since $R$ is $m$-perfect). Then,
$\Ext(M,H)=0$. Thus, we have the following exact sequence:
$$0 \longrightarrow \Hom(M,H) \longrightarrow \Hom(M,P^0)
\stackrel{\Hom(M,f)}
 \longrightarrow \Hom(M,F) \longrightarrow \Ext(M,H)=0.$$
Then, there exists $\overline{\varphi}:\, M \rightarrow P^0$ such
that $\varphi=f\overline{\varphi}$. Since $\varphi$ is injective,
$\overline{\varphi}$ is also injective, and so we obtain the
following short exact sequence of $R$-modules: $$(*)\qquad 0
\rightarrow M \stackrel{\overline{\varphi}}\rightarrow P^0
\rightarrow G^0 \rightarrow 0.$$ \indent Now, to complete the
proof, it remains to prove that $\Ext^{i}(G^0,F')=0$ for all $i>0$
and all $R$-modules $F'$ with finite projective dimension.\\
\indent First, assume that $F'$ is projective. Since $\varphi $ is
a flat preenvelope of $M$, there exists, for all $\alpha \in
\Hom(M,F')$, a homomorphism $g:\, F\rightarrow F'$ such that
$\alpha=g\varphi $, hence $\alpha=g f\overline{\varphi}$. This
means that the functor $\Hom(-,F')$ leaves the short  sequence
$(*)$ exact. Then, by the long exact sequence $$0 \rightarrow
\Hom(G^0,F') \rightarrow \Hom(P^0,F') \rightarrow \Hom(M,F')
\rightarrow \Ext(G^0,F')\rightarrow  \Ext( P^0,F')=0,$$ we deduce
that $\Ext(G^0,F')=0$. Also, we use the short exact sequence $(*)$
to deduce that $\Ext^{i}(G^0,F')=0$ for all $i>0$ and all
projective $R$-modules $F'$. Finally, this implies directly that
$\Ext^{i}(G^0,F')=0$ for all $i>0$ and all $R$-modules $F'$ with
finite projective dimension.\cqfd\bigskip

\noindent\textbf{Proof of Theorem \ref{thm-main}.} First, from
\cite[Theorem 7.2.5 (2)]{DM} and \cite[Theorem 2.28]{HH}, the
inequality $\cot.D(R)\leq\Ggldim(R)$ holds for any arbitrary ring
$R$.\\
\indent Then, we prove the inequality
$\Ggldim(R)\leq\Gwdim(R)+\cot.D(R)$ when $R$ is coherent. For
that, we may assume that $\cot.D(R)=m$ and $\Gwdim(R)=n$ are
finite (i.e., $R$ is $m$-perfect and $n$-FC). Let $M$ be an
$R$-module, and consider an exact sequence of $R$-modules:
$$0\rightarrow K_n \rightarrow P_{n-1}\rightarrow\cdots
\rightarrow P_0\rightarrow M \rightarrow 0,$$ where each $P_i$ is
projective, and, from \cite[Theorem 4.11]{BM1}, $K_n$ is
Gorenstein flat. We have:
$$(*)\quad \Ext^{n+k}(M,Q)\cong\Ext^{k}(K_n,Q)\qquad  \mathrm{for}\ \mathrm{all}
\ k\geq 1\ \mathrm{and}\ \mathrm{all}\ \mathrm{modules}\ Q.$$
Assume that $Q$ is a projective $R$-module, and consider an exact
sequence of $R$-modules:
$$0\rightarrow Q \rightarrow C_0\rightarrow \cdots \rightarrow
C_{m-1}\rightarrow C_{m} \rightarrow 0,$$ where $C_i$ is injective
for $i=1,...,m-1$, and then, from \cite[Proposition 7.2.1]{DM},
$C_m$ is cotorsion. We have:
$$(**)\quad \Ext^{m+i}(K_n,Q)\cong \Ext^{i}(K_n,C_m)\qquad
\mathrm{for}\ \mathrm{all}\ i\geq1.$$ Since $\Gwdim(R)$ is finite,
each of the $R$-modules $C_0$,...,$C_{m-1}$ has finite flat
dimension (from \cite[Theorem 4.11]{BM1}). Then, $C_m$ has finite
flat dimension. Thus, $\Ext^{i}(K_n,C_m)=0$ for all $i\geq 1$
(from \cite[Proposition 3.22]{HH} and since  $K_n$ is Gorenstein
flat). Then, by $(*)$ and $(**)$, $\Ext^{n+m+i}(M,Q)=0$ for all
$i\geq1$. This implies, from Lemma \ref{lem-Gpro-n-perfect-FC},
that $\Gpd(M)\leq n+m$, as desired.\cqfd\bigskip

Theorem \ref{thm-main} enables us to compute the Gorenstein global
dimension of some particular cases of trivial extensions of rings
and of group rings.\bigskip

Recall that the trivial extension of a ring $R$ by an $R$-module
$M$ is the ring denoted by $R\ltimes M$ whose underling group is
$A\times M$ with multiplication given by
$(r,m)(r',m')=(rr',rm'+r'm)$ (see for instance \cite{Triv} and
\cite[Chapter 4, Section 4]{Glaz}). Next result compute the
Gorenstein global dimension of a particular case of trivial
extensions of rings. For that, we use the notion of finitistic
projective dimension of rings. Recall the finitistic projective
dimension of a ring $R$, denoted by $\FPD(R)$, is defined by:
$$\FPD(R)=\sup\{\pd_R(M)|M\ R\!-\!\mathrm{module}\ \mathrm{with}\
\pd_R(M)<\infty\}.$$ From  \cite[Theorem 2.28]{HH}, we have for
every ring $R$: $\FPD(R)\leq \Ggldim(R)$, with equality if
$\Ggldim(R)$ is finite.

\begin{prop}\label{pro-Ggldim-trivial} Let $R \ltimes
R$ be the trivial extension of a ring  $R$ by  $R$. Then,
$\FPD(R\ltimes R)= \FPD(R)$, $\cot.D(R\ltimes R)= \cot.D(R)$,
and $\gldim(R\ltimes R)=\infty$.\\
\indent Furthermore, if $R$ is coherent, then  $\Ggldim(R\ltimes
R)= \Ggldim(R)$.
\end{prop}

The proof of this theorem involves the following results:

\begin{lem}[\cite{Triv},  Theorem 4.28 and Remark page 81]\label{lem-triv1}
Let $R$ be a ring and let $M$ be any non-zero cyclic $R$-module.
Then, $\FPD(R\ltimes M)= \sup\{\pd_R(N)<\infty:\ \Tor_i^R(M,N)=0\
\mathrm{for}\ \mathrm{all}\ i>0\}$ and $\gldim(R\ltimes
M)=\infty$.\\
\indent In particular, if $F$ is an $R\ltimes M$-module having
finite projective dimension, then $\pd_{R\ltimes
M}(F)=\pd_R(R\otimes_{R\ltimes M} F)$.
\end{lem}

\begin{lem}[\cite{Triv},  Theorem  4.32]\label{lem-triv2}
Let $R$ be a ring and let $M$ be an $R$-module such that:

$$\Ext^i_R(M,M)\cong \left\{%
\begin{array}{ll}
    R & \hbox{if\ i=0;} \\
   0 & \hbox{if\ i$>$0.} \\
\end{array}%
\right.$$ Then, $\id_{R\ltimes M}(R\ltimes M)=\id_{R }( M)$.
\end{lem}

\noindent\textbf{Proof of Proposition \ref{pro-Ggldim-trivial}.}
From Lemma \ref{lem-triv1},
$\FPD(R\ltimes R)= \FPD(R)$ and $\gldim(R\ltimes R)=\infty$.\\
\indent  We prove that $\cot.D(R\ltimes R)= \cot.D(R)$. From Lemma
\ref{lem-triv1}, we have:
$$\pd_{(R\ltimes R)}(F)=\pd_{R}(R\otimes_{(R\ltimes R)} F )$$ for
every $R\ltimes R$-module $F$ with finite projective dimension.
This implies that $\cot.D(R\ltimes R)\leq \cot.D(R)$. Conversely,
consider a flat $R$-module $F$, then $F \otimes_{R} (R\ltimes R) $
is a flat $R\ltimes R$-module. Thus, since $ R\ltimes R$ is a free
$R$-module such that $ R\ltimes R \cong_R R^2$ we have:
$$\pd_{R}(F\otimes_{R} R)= \pd_{R}(F\otimes_{R} (R\ltimes R))\leq
\pd_{(R\ltimes R)}(F\otimes_{R} (R\ltimes R)) \leq \cot.D(R\ltimes
R).$$ Therefore, $ \cot.D(R)\leq \cot.D(R\ltimes R)$, as
desired.\\ \indent Now, we assume that $R$ is coherent, and we
prove the equality $\Ggldim(R\ltimes R)= \Ggldim(R)$.\\ \indent
First assume that $\Ggldim(R)$ is finite. Then, by the reason
above and from \cite[Theorem 2.28]{HH}, $\FPD(R\ltimes
R)=\FPD(R)=\Ggldim(R)$ is finite. So $\cot.D(R\ltimes R)$ is
finite. On the other hand, from
 Lemma \ref{lem-triv2}, $\id_{(R\ltimes R)}(R\ltimes
R)=\id_R(R) $ which is finite (by \cite[Lemma 3.3]{BM1} and since
$\Ggldim(R)$ is finite). Then, by \cite[Theorem 4.11]{BM1},
$\Gwdim(R)= \FPid_R(R)\leq \id_R(R)$ is finite. Then, from Theorem
\ref{thm-main}, $\Ggldim(R\ltimes R)$ is finite. Therefore,
from \cite[Theorem 2.28]{HH}, $\Ggldim(R\ltimes R)= \FPD(R\ltimes R)= \Ggldim(R)$.\\
\indent Similarly we show that  $\Ggldim(R\ltimes R)= \Ggldim(R)$
when $\Ggldim(R\ltimes R)$  is finite, and this gives the desired
result.\cqfd\bigskip

As mentioned in the introduction,  the Noetherian rings of finite
Gorenstein global dimension are the same Iwanaga-Gorenstein rings;
and in the class of rings of finite weak dimension the global
dimension and the Gorenstein global dimension coincide. In the
following example, we construct a family of non-Noetherian
coherent rings $\{S_i\}_{i\geq 1}$ such that $\Ggldim(S_i)=i$ and
$\wdim(S_i)=\infty$ for every $i\geq 1$.

\begin{exmp}\label{exm-3}
Let $R_n=R[X_1,X_2,...,X_n]$ be the polynomial ring in $n$
indeterminates over  a non-Noetherian hereditary ring $R$. Let
$S_i=  R_{i-1}  \ltimes R_{i-1}$ be the trivial extension of
$R_{i-1}$ by $R_{i-1}$ for $i\geq 1$ (such that $R_0=R$). Then,
for every $i\geq 1$, $S_i $ is a non-Noetherian coherent ring with
$\Ggldim(S_i )=i$ and $\wdim(S_i)=\infty$.
\end{exmp}
\proof  From \cite[Theorem 7.3.1]{Glaz}, $R_n=R[X_1,X_2,...,X_n]$
is coherent for every $n\geq 1$. And by Hilbert's syzygy theorem,
$\gldim(R_n)=\gldim(R)+n=1+n$. Therefore, Proposition
\ref{pro-Ggldim-trivial} implies that  $\Ggldim(S_i )=i$ for every
$i\geq 1$.\\
\indent Finally,   $\wdim( S_i)=\infty$ for every $i\geq 1$
follows from \cite[Proposition 3.11]{BM1} and since $\gldim(S_i
)=\infty$ from Proposition \ref{pro-Ggldim-trivial}.\cqfd\bigskip


We end this paper with a study of the global cotorsion dimension
of group rings, and then the Gorenstein global dimension of a
particular group ring is computed.\\
\indent Let $R$ be a ring and let $G$ be an abelian group written
multiplicatively. The free $R$-module on the elements of $G$ with
multiplication induced by $G$ is a ring, called group ring of $G$
over $R$ and denoted by $RG$ (see for instance \cite[Chapter 8, Section 2]{Glaz}).\\
\indent In \cite{Perf-RG}, we have that $RG$ is perfect if and
only if $R$ is perfect and $G$ is finite. Here, we set the
following extension.

\begin{thm}\label{thm-RG-Cot}
Let $R$  be a ring and let $G$ be an abelian group. We have:
$$\cot.D(R)\leq \cot.D(RG)\leq\cot.D(R)+\pd_{RG}(R).$$
\indent Furthermore, if $G$ and $\pd_{RG}(R)$ are finite, then
$\cot.D(R)= \cot.D(RG)$.
\end{thm}

To prove this result we need the following lemma.

\begin{lem} [\cite{CE}, page 352]\label{lem-RG-Cot}
Let $R$  be a ring, let $G$ be an abelian group, and let $M$ and
$N$ be two $RG$-modules satisfying $\Ext^{p}_{R}(M,N)=0$ for all
$p>0$. Then, $$\Ext^{n}_{RG}(M,N)\cong
\Ext^{n}_{RG}(R,\Hom_R(M,N))$$ for all $n>0$, where
$\Hom_{R}(M,N)$ is the $RG$-module defined by
$(gf)(x)=g[f(g^{-1}x)]$  for $x\in M$, $f\in\Hom_{R}(M,N)$, and
$g\in G$.
\end{lem}

\noindent\textbf{Proof of Theorem \ref{thm-RG-Cot}.} First, we
prove the inequality  $\cot.D(R)\leq \cot.D(RG)$. We may assume
that $\cot.D(RG)=n$ is finite. Let $F$ be a flat $R$-module, then
$F\otimes_R RG$ is a flat $RG$-module.  Since $RG\cong R^{(G)}$ is
a free $R$-module, $\pd_{R}(F)=\pd_{R}(F^{(G)})=\pd_{R}(F\otimes_R
RG)\leq \pd_{RG}(F\otimes_R RG)\leq n.$ This implies the desired
inequality.\\ \indent Now, we prove the inequality
$\cot.D(RG)\leq\cot.D(R)+\pd_{RG}(R)$. For that we may assume that
$\cot.D(R)=s$ and  $\pd_{RG}(R)=r$ are finite. Let $F$ be a flat
$RG$-module (then it is also flat as an $R$-module), and consider
an exact sequence of $RG$-modules: $$0\rightarrow P_s\rightarrow
\cdots \rightarrow P_0\rightarrow F \rightarrow 0,$$ where
$P_0$,...,$P_{s-1}$ are projective $RG$-modules, then they are
projective as $R$-modules, and so $P_s$ is a projective $R$-module
(since $\cot.D(R)=s$). Thus, $\Ext^{p}_{R}(P_s,N)=0$ for all $p>0$
and all $R$-modules $N$. Then, from Lemma \ref{lem-RG-Cot} above
and since $\pd_{RG}(R)=r$, $$\Ext^{n}_{RG}(P_s,N)\cong
\Ext^{n}_{RG}(R,\Hom_R(P_s,N))=0$$ for all $n>r$ and all
$RG$-modules $N$. Thus, $\pd_{RG}(P_s)\leq r$ and so
$\pd_{RG}(F)\leq s+r$. Therefore, $\cot.D(RG)\leq s+r$, as
desired.\\ \indent Assume now that  $G$   and $\pd_{RG}(R)$ are
finite. From \cite[Lemma 3.2 (a)]{BG}, $R$ is projective as an
$RG$-module, and by the inequalities above,  $\cot.D(R)=
\cot.D(RG)$, as desired.\cqfd\bigskip

The above result  and the main result (Theorem \ref{thm-main}) are
used to compute the Gorenstein global dimension of a particular
group ring as follows:

\begin{prop}\label{prop-RG-appl} Let $R$ be a ring with $\Gwdim(R)=0$.
If $G$ is a finite group such that its order is invertible in $R$,
then $\Gwdim(RG)=0$ and $\Ggldim(RG)=\Ggldim(R)$.
\end{prop}
\proof  First, note that $R$ is coherent and so it is an IF-ring
(from \cite[Theorem 6]{FCring} and \cite[Proposition
4.2]{Stenst}). Then, from \cite[Theorem 3  page  250]{Colby}, $RG$
is
an IF-ring and so $\Gwdim(RG)=0$.\\
\indent Now, by Theorem \ref{thm-main}, $\Ggldim(R) = \cot.D(R)$
and $\Ggldim(RG) = \cot.D(RG)$. And from \cite[Theorem
8.2.7]{Glaz}, $R$ is projective as $RG$-module. Thus, from Theorem
\ref{thm-RG-Cot}, $\cot.D(RG)= \cot.D(R)$. This implies the
desired equality $\Ggldim(RG)=\Ggldim(R)$.\cqfd\bigskip

\end{section}

\end{document}